\documentclass[12pt]{amsart}
\usepackage{amssymb,latexsym,esint}
\usepackage{enumerate}
\usepackage{float}
\usepackage{import}

\usepackage{mathtools}
\usepackage{amsfonts}
\usepackage{amssymb}
\usepackage{amsmath}
\usepackage{graphicx}
\usepackage{graphicx,color}
\usepackage[normalem]{ulem}
\usepackage{todonotes}
\usepackage{mathrsfs}
\usepackage{import}
\usepackage{xifthen}
\usepackage{pdfpages}
\usepackage{transparent}
\usepackage{yfonts}

\makeatletter
\@namedef{subjclassname@2020}{%
  \textup{2020} Mathematics Subject Classification}
\makeatother

\usepackage{amsmath,amstext,amsthm,amsxtra,mathtools,amssymb}

\usepackage[bookmarksopen,bookmarksdepth=3,colorlinks,citecolor=red,pagebackref,hypertexnames=true]{hyperref}
\usepackage[backrefs,msc-links,nobysame,initials,non-sorted-cites]{amsrefs}

\usepackage[normalem]{ulem}
\usepackage{lmodern}

\usepackage{lmodern}

\usepackage[backrefs,msc-links,nobysame,initials]{amsrefs}

\allowdisplaybreaks

\newtheorem{thm}{Theorem}

\newtheorem{lem}{Lemma}

\theoremstyle{definition}

\theoremstyle{remark}

\begin{document}

\title{Sharp Weak Type Estimates for Maximal Operators associated to Rare Bases}
\dedicatory{Dedicated to Dmitriy Dmitrishin on the occasion of his fifty-fifth birthday}

\author{Paul Hagelstein}
\address{P. H.: Department of Mathematics, Baylor University, Waco, Texas 76798}
\email{\href{mailto:paul_hagelstein@baylor.edu}{paul\_hagelstein@baylor.edu}}
\thanks{P. H. is partially supported by a grant from the Simons Foundation (\#521719 to Paul Hagelstein).}

\author{Giorgi Oniani}
\address{G. O.: School of Computer Science and Mathematics, Kutaisi International University, Youth Avenue, Turn 5/7, Kutaisi 4600, Georgia}
\email{\href{mailto:giorgi.oniani@kiu.edu.ge}{giorgi.oniani@kiu.edu.ge}}

\author{Alex Stokolos}
\address{A. S.: Department of Mathematical Sciences, Georgia Southern University, Statesboro, Georgia 30460}
\email{\href{mailto:astokolos@GeorgiaSouthern.edu}{astokolos@GeorgiaSouthern.edu}}

\subjclass[2020]{Primary 42B25}
\keywords{maximal operator, weak type estimate, basis}

\begin{abstract} Let $\mathcal{B}$ denote a nonempty translation invariant collection of intervals in $\mathbb{R}^n$ (which we regard as a rare basis), and define the associated geometric maximal operator $M_\mathcal{B}$ by
$$M_\mathcal{B}f(x) = \sup_{x \in R \in \mathcal{B}} \frac{1}{|R|}\int_R |f|.$$
We provide a sufficient condition on $\mathcal{B}$ so that 
the estimate
$$	|\{x \in \mathbb{R}^n : M_{\mathcal{B}}f(x) > \alpha\}|\leq C_n \int_{\mathbb{R}^{n}} \frac{|f|}{\alpha}\left(1+\log^+\frac{|f|}{\alpha}\right)^{n-1}
$$
is sharp. As a corollary we obtain sharp weak type estimates for maximal operators \mbox{associated} to several  classes of rare bases including  C\'ordoba, Soria  and Zygmund bases. 

\end{abstract}

\maketitle
\section{Introduction}

Let $\mathcal{B}$ be a \emph{basis} in $\mathbb{R}^n$, which we can treat as a collection of sets of positive finite measure covering $\mathbb{R}^n$.  We may associate to $\mathcal{B}$ a \emph{geometric maximal operator} $M_\mathcal{B}$ defined on measurable functions $f$ on $\mathbb{R}^n$ by
$$M_\mathcal{B}f(x) = \sup_{x \in R \in \mathcal{B}} \frac{1}{|R|}\int_R |f|\;.$$
Two well-known examples of geometric maximal operators include the \emph{Hardy-Littlewood \mbox{maximal} operator} $M_{HL}$ and the \emph{strong maximal operator} $M_{S}$.  For the Hardy-Littlewood maximal operator, the basis $\mathcal{B}$  consists of all cubic intervals in $\mathbb{R}^n$; for the strong maximal operator, the basis consists of all intervals in $\mathbb{R}^n$. (For clarity, in this paper an \mbox{\emph{interval}} in $\mathbb{R}^n$ is a rectangular parallelepiped whose sides are parallel to the coordinate axes.)

The Hardy-Littlewood maximal operator satisfies the weak type $(1,1)$ estimate
$$|\{x \in \mathbb{R}^n : M_{HL}f(x) > \alpha\}| \leq C_n\int_{\mathbb{R}^n}\frac{|f|}{\alpha}.$$  The strong maximal operator satisfies the weaker  estimate
\begin{equation}\label{e1}|\{x \in \mathbb{R}^n : M_{S}f(x) > \alpha\}| \leq C_n \int_{\mathbb{R}^n}\frac{|f|}{\alpha} \left(1 + \log^+ \frac{|f|}{\alpha}\right)^{n-1}\end{equation}
and moreover this inequality is \emph{sharp} in the sense that, if $\phi:[0, \infty) \rightarrow [0,\infty)$ is a convex increasing function with $\lim_{u \rightarrow \infty}\frac{\phi(u)}{u (1 + \log u)^{n-1}} = 0$, then there is no finite constant $C_{n, \phi}$ such that the estimate $$|\{x \in \mathbb{R}^n : M_{S}f(x) > \alpha\}| \leq C_{n, \phi} \int_{\mathbb{R}^n}\phi \left(\frac{|f|}{\alpha}\right)$$ holds for all measurable functions $f$ and $\alpha > 0$.
\\

Geometric maximal operators associated to \emph{rare bases} occupy a fascinating middle ground between the Hardy-Littlewood and strong maximal operators.  Significant  mathematical work on the topic on rare bases has been done by, among others, Zygmund \cite{zygmund1967}, C\'ordoba \cite{cordoba}, \mbox{Soria \cite{soria},} and  Rey \cite{rey}. For the purposes of this paper, a rare basis is a translation invariant collection of some (but not  all) intervals in $\mathbb{R}^n$.  The problem of finding a sharp weak type $\phi(L)$ estimate for the maximal operator $M_\mathcal{B}$ associated to  a rare basis $\mathcal{B}$ is of central importence in the area. A natural conjecture (see, e.g., \cite{hs2022pams}) is  $M_\mathcal{B}$ must satisfy a sharp weak type $L(1 + \log^{+} L)^{k}$ estimate for some  $k\in \{0,1,\dots, n-1\}$.  This conjecture in the two-dimensional case was proven by Stokolos \cite{stokolos1988} (see also \cite{stokolos2005}).  Note that for any rare basis $\mathcal{B}$, the associated geometric maximal operator $M_\mathcal{B}$ is dominated by the strong maximal operator $M_S$ and hence $M_\mathcal{B}$ automatically satisfies the  weak type $L(1 + \log^{+} L)^{n-1}$ estimate (\ref{e1}). Results of Stokolos \cite{stokolos1988, stokolos2006};  Dmitrishin, Hagelstein, and Stokolos \cite{dhs2022, hs2022pams, hs2022cordoba} and  D'Aniello and Moonens \cite{dm2017} have provided examples of a variety of multi-dimensional rare bases in which this worst-case estimate is sharp.  The purpose of this paper is to present a generalization of results in these papers that provides a set of conditions on a rare basis $\mathcal{B}$ that will guarantee the sharpness of the weak type  $L(1 + \log^{+} L)^{n-1}$ estimate on $M_\mathcal{B}$.

In Section 2 we will state the main theorem of this paper and provide the associated requisite terminology.   Section 3 will be devoted to a proof of the theorem.  In \mbox{Section 4} we will provide a useful generalization of 
our main theorem and subsequently provide applications of our results, in particular   showing how they imply that C\'ordoba, Soria and Zygmund bases in $\mathbb{R}^n$  generate geometric maximal operators for which the  weak type $L(1 + \log^{+}L)^{n-1}$ estimate  is sharp.

\section{Terminology and Statement of Main Theorem}

Let $\mathcal{B}$ be a rare basis in $\mathbb{R}^n$. The \textit{spectrum} of   $\mathcal{B}$  will be defined as the set of all $n$-tuples of the type
$$
(\lceil \log |R_1| \rceil, \dots, \lceil \log |R_n| \rceil)
$$
where $R_1 \times \dots \times R_n \in \mathcal{B}$, $\lceil x  \rceil$ denotes the least integer greater than or equal to $x$ and here and below $\log$ stands for $\log_2$.  The spectrum will be denoted by $W_\mathcal{B}$.

Let us call a set   $W\subset \mathbb{Z}$  a  \textit{net} for a  set  $S\subset \mathbb{Z}$  if there exists a number $N \in \mathbb{N}$ such that for every  $s\in S$ there exists $w\in W$ with $|s-w|\leq N$.

Given  a set  $W\subset \mathbb{Z}^n$ and $t\in \mathbb{Z}^{n-1}$, we let $W_t$ denote the set $\{\tau\in \mathbb{Z}: (t,\tau)\in W \}$.

Let us call a set  $W\subset \mathbb{Z}^n$  a \textit{net} for  a set $S\subset \mathbb{Z}^n$  if $W_t$ is a net for $S_t$ for every $t\in  \mathbb{Z}^{n-1}$. 

If $k \in \{1, 2, \dots, n\}$, we let $\pi_k: \mathbb{R}^n \rightarrow \mathbb{R}^k$  denote the projection defined by $\pi_k(x_1,\dots, x_n) =  (x_{1},\dots, x_{k})$ .  

If $S_1,\dots, S_n\subset \mathbb{Z}$, we say that the set $W\subset \mathbb{Z}^n$ is  \textit{dense} in the set $S_1\times \dots \times S_n$ provided
the sets 
$$
\pi_{1}(W), \; \pi_{2}(W), \; \dots, \; \pi_{n}(W)=W
$$
are nets for the sets
$$
S_1, \; \pi_{1}(W)\times S_2,\; ..., \; \pi_{n-1}(W)\times S_n,
$$
respectively.

\mbox{ }

Our main theorem is the following.

\medskip

\begin{thm}\label{t1} If for a rare basis  $\mathcal{B}$ in $\mathbb{R}^n$ there exist  infinite sets $S_1,\dots, S_n \subset \mathbb{Z}$ for which the	 spectrum of $\mathcal{B}$ is dense in $S_1\times \dots \times S_n$,  then  the maximal operator $M_\mathcal{B}$ satisfies a sharp weak type $L(1+ \log^{+}L)^{n-1}$  estimate. Moreover, 
	for every $\alpha \in(0,1)$ there exists a bounded set $E_{\alpha}\subset \mathbb{R}^n$ with  positive  measure such that
$$
|\{x \in \mathbb{R}^n : M_{\mathcal{B}} (\chi_{E_{\alpha}})(x) > \alpha\}|\geq c_{n} \frac{1}{\alpha} \bigg(1 + \log \frac{1}{\alpha}\bigg)^{n-1}|E_{\alpha}|.
$$

\end{thm}


\section{Proof of Theorem \ref{t1}}

 \begin{lem}\label{l1} Let $\mathcal{B}$ be a rare basis  and suppose each  of the intervals from $\mathcal{B}$ has dyadic side lengths. Suppose $S_1,\dots, S_n$ $\subset \mathbb{Z}$ are infinite sets and the	 spectrum of $\mathcal{B}$ is dense in \mbox{$S_1\times \dots \times S_n$.} Then for every $k\in \mathbb{N}$ there exist increasing sequences $(s_{1,m})_{m=0}^k, \dots, (s_{n,m})_{m=0}^k$ with members from $S_1,\dots, S_n$ respectively such that for every $n$-tuple  $(m_1,\dots, m_n)$ \mbox{belonging} to $\{1,\dots,k\}^n$ there exists an interval $R_1\times \dots \times R_n$ from the basis $\mathcal{B}$ for which $|R_1|\in(2^{s_{1, m_1-1}}, 2^{s_{1, m_1}}],\dots,|R_n|\in(2^{s_{n, m_n-1}}, 2^{s_{n, m_n}}].
$
\end{lem}

\begin{proof}  Since  $\pi_{1}(W_\mathcal{B})$ is a net for $S_1$, there exists a number $N$ such that for every  $s\in S_1$ there exists $\tau\in \pi_{1}(W_\mathcal{B})$ with $|s-\tau|\leq N$. Let $\alpha_{1,0}< \dots < \alpha_{1,2k}$ be numbers from $S_1$ such that $\alpha_{1,m}-\alpha_{1,m-1}>N$ for every $m\in\{1,\dots, 2k\}$. Set $s_{1,m}=\alpha_{1,2m}$ ($m\in \{0,\dots, k\}$). Then it is easy to see that  
\begin{equation}\label{e2}
\pi_{1}(W_\mathcal{B})\cap (s_{1,m-1},s_{1,m}]\neq \emptyset 
\end{equation}
for every $m\in \{1,\dots,k\}$.

Suppose that for some $j<n$ increasing sequences $(s_{1,m})_{m=0}^k, \dots,$ $(s_{j,m})_{m=0}^k$  with members from $S_1,\dots, S_j$ respectively are constructed.
  
Let us consider an arbitrary $(t_1,\dots,t_j)\in \pi_{j}(W_\mathcal{B})$ such that $t_1\in [s_{1,0}, s_{1,k}],$ $\dots, \mbox{$t_j\in [s_{j,0}, s_{j,k}]$}$.
Let
$$
W_{\mathcal{B}, t_1,\dots,t_j}=\{\tau\in \mathbb{Z}: (t_1,\dots, t_j,\tau)\in \pi_{j+1}(W_\mathcal{B}) \}.
$$
Since  $W_\mathcal{B}$ is dense in  $S_1\times \dots \times S_n$ we have $W_{\mathcal{B},t_1,\dots,t_j}$ is a net for $S_{j+1}$. Let $N_{t_1,\dots, t_j}$ be the number such that for every $s\in S_{j+1}$ there exists $\tau\in W_{\mathcal{B},t_1,\dots,t_j} $  for which $|s-\tau|\leq N_{t_1,\dots, t_j}$.

Denote by $N$ the largest of the numbers $N_{t_1,\dots, t_j}$ where $(t_1,\dots,t_j)$ $\in \pi_{j}(W_\mathcal{B})$ and $t_1\in [s_{1,0}, s_{1,k}],$ $\dots, t_j\in [s_{j,0}, s_{j,k}]$.  Let $\alpha_{j+1,0}< \dots < \alpha_{j+1,2k}$ be numbers from $S_{j+1}$ such that $\alpha_{j+1,m}-\alpha_{j+1,m-1}>N$ for every $m\in\{1,\dots, 2k\}$. Set $s_{j+1,m}=\alpha_{j+1,2m}$ ($m\in \{0,\dots, k\}$). Then  for every $(t_1,\dots,t_j)$ $\in \pi_{j}(W_\mathcal{B})$ with $t_1\in [s_{1,0}, s_{1,k}], \dots, t_j\in [s_{j,0}, s_{j,k}]$ we have that  
\begin{equation}\label{e3}
W_{\mathcal{B},t_1,\dots,t_j}\cap (s_{j+1,m-1},s_{j+1,m}]\neq \emptyset 
\end{equation}
for every $m\in \{1,\dots,k\}$.

Taking into account (\ref{e2}) and (\ref{e3}), it is easy to check that the sequences $(s_{1,m})_{m=0}^k, \dots,$ $(s_{n,m})_{m=0}^k$ constructed in such a way   have the  needed property. The lemma is proved. 
\end{proof}

\medskip

Suppose $E\subset \mathbb{R}$ is a measurable set, $\Omega$ is a collection of disjoint closed one-dimensional intervals, $H=\bigcup_{I\in\Omega} I$ and $0\leq \alpha\leq 1$.  We will say that $E$ $\alpha$-\textit{saturates} $H$ (notation: $E\leftarrow^{\alpha} H$) if $|I\cap E|/|I|=\alpha$ for every $I\in\Omega$.  
 For the case $\alpha=1/2$ we will write simply $E\leftarrow H$. 
 
 Note that: 
 
 a) If $E\leftarrow^{\alpha} H$ then $|H\cap E|/|H|=\alpha$;
 
 b) If $E\subset H$, $E\leftarrow^{\alpha} H$ and $H\leftarrow^{\beta} T$ then $E\leftarrow^{\alpha\beta} T$.

Let $I$ and $J$ be closed one-dimensional intervals with lengths $2^{p}$ and $2^{q}$ where $p, q\in \mathbb Z$ and $p<q$. Let us consider the partition of $J$ into  non-overlapping closed subintervals $J_1,\dots,J_{2^{p-q}}$ having the same length as $I$ and such that $\min J_1<\dots < \min J_{2^{p-q}}$. By $\langle I, J \rangle$ denote the union of the intervals $J_k$ with odd indices. Obviously, $\langle I, J \rangle \subset J$ and $\langle I, J \rangle \leftarrow J$.

Let $I$ be a closed one-dimensional interval with a length $2^{p}$  where $p\in \mathbb Z$,  $\Omega$ be a collection of disjoint closed one-dimensional intervals having a length $2^{q}$ with $q\in \mathbb Z$  and $p < q$, and  $H=\bigcup_{J\in\Omega} J$. By $\langle I, H \rangle$ we denote the union  $\bigcup_{J\in\Omega}\langle I, J \rangle$.  Obviously, $\langle I, H \rangle \subset H$ and \mbox{ $\langle I, H \rangle \leftarrow H$. }

It is easy to check the validity of the following statement.
\medskip

\begin{lem}\label{l2} \textit{Suppose $s_0<s_1<\dots<s_k$ are some  integers, $I_0=$ $[0,2^{s_0}],\dots, I_k=[0,2^{s_k}]$ and} 
$$
I^{\ast}_k=I_k,\; I^{\ast}_{k-1} = \langle I_{k-1}, I^{\ast}_k \rangle, \; I^{\ast}_{k-2} = \langle I_{k-2},I^{\ast}_{k-1}\rangle, \;\dots,\; I^{\ast}_{0} =  \langle I_0, I^{\ast}_{1} \rangle. 
$$
\textit{Then the sets $I^{\ast}_0,  I^{\ast}_1, \dots,  I^{\ast}_k$ have the following properties}: 

1)  $I^{\ast}_0 \subset  I^{\ast}_1 \subset \dots \subset  I^{\ast}_k$;

2) \textit{$I^{\ast}_0 \leftarrow  I^{\ast}_1 \leftarrow \dots \leftarrow  I^{\ast}_k $ and, moreover, $I^{\ast}_{m-1} \leftarrow J\; $ for every $ m \in \{1,\dots, k\}$ and  every dyadic interval $J$ contained in $I^{\ast}_m$ whose length belongs to  $(2^{s_{m-1}},2^{s_{m}}]$};

3) \textit{$I^{\ast}_0 \leftarrow^{1/2^m}  I^{\ast}_m\; $ for every $m\in\{1,\dots, k\}$ and, moreover,  $I^{\ast}_0 \leftarrow^{1/2^m}  J\; $ for every \mbox{$m\in\{1,\dots, k\}$} and every dyadic interval $J$ contained in $I^{\ast}_m$ whose length belongs to  $(2^{s_{m-1}}, 2^{s_{m}}]$. }
\end{lem}

\medskip


 \noindent \textbf{Remark 1}. It is important to recognize that the interval $J$ in 3) above may have not only  length $2^{s_m}$ but also the ``intermediate'' lengths  $2^{s_{m-1}+1}, \dots,   2^{s_{m}-1}$ as well. 
\\

\noindent \textbf{Remark 2}. The idea of using the sets  $I^{\ast}_0,  I^{\ast}_1, \dots,  I^{\ast}_k$ in Lemma \ref{l2} for the study of weak type estimates for maximal operators associated to rare bases goes back to \mbox{Stokolos  \cite{stokolos2006}. }

\medskip

For every $k\geq n$, by $\Omega_{n,k}$ we denote the set of all $n$-tuples  $ (m_1,\dots,$ $  m_n)\in \mathbb{N}^n$ for which  $m_1+\dots+$ $m_n=k$. Clearly, $\textup{card}(\Omega_{n,k})\leq k^{n-1}$. On the other hand, it is easy to see that $\textup{card}(\Omega_{n,k})\geq c_n k^{n-1}$.

\medskip


\begin{lem}\label{l3}
For every $k \geq n$ and increasing sequences of integers $$
	(s_{1,m})_{m=0}^k,\;\dots\;,(s_{n,m})_{m=0}^k
	$$
		there exists a bounded set $E \subset \mathbb{R}^n$ of positive measure with the following property:   if   $\mathcal{B}$ is a rare basis in $\mathbb{R}^n$  and $\Omega$ is the set of all $n$-tuples  $ (m_1,\dots,$ $  m_n) \in \Omega_{n,k}$ for which  there exists an interval $R_1\times \dots \times R_n \in B$ with dyadic side lengths such that 
	$$
	|R_1|\in(2^{s_{1, m_1-1}}, 2^{s_{1, m_1}}],\dots,|R_n|\in(2^{s_{n, m_n-1}}, 2^{s_{n, m_n}}],
	$$
then
$$
	|\{x \in \mathbb{R}^n : M_{\mathcal{B}} (\chi_{E})(x)\geq 1/2^k\}|\geq c_{n}\, \textup{card}(\Omega) \, 2^{k}|E|.
	$$
\end{lem}
\mbox{ }

\begin{proof} For every $j\in \{1,\dots, n\}$ denote  by  $I_{j,0}, \dots,I_{j,k}$ and $I^{\ast}_{j,0}, \dots,I^{\ast}_{j,k} $   the sets  corresponding to the   sequence $s_{j,0},\dots, s_{j,k}$ according to  Lemma 2. We define the set $E \subset \mathbb{R}^n$  by
$$E=I^\ast_{1,0}\times\dots\times I^\ast_{n,0}\;.
$$

We now show that 
\begin{equation}\label{e4}
\{x \in \mathbb{R}^n : M_{\mathcal{B}} (\chi_{E})(x)\geq 1/2^k\}\supset \bigcup_{(m_1,\dots, m_n)\in \Omega}  
I^{\ast}_{1,m_1}\times \dots \times I^{\ast}_{n,m_n}. 
\end{equation}

Let us consider an arbitrary  $(m_1,\dots, m_n) \in\Omega$. Let $R_1\times \dots \times R_n$ be an interval from $\mathcal{B}$ such that  $|R_1|,\dots, |R_n|$ are dyadic numbers and   $|R_j|\in(2^{s_{j, m_j-1}}, 2^{s_{j, m_j}}]$ for every $j\in \{1,\dots, n\}$. Then for every $j\in\{1,\dots, n\}$  each component interval of $I^{\ast}_{j,m_j}$ can be decomposed into pairwise non-overlapping subintervals $\Delta$ each of whose lengths is equal to $|R_j|$. By Lemma \ref{l2}  (see  statement 3)) each $\Delta$ will be $1/2^{m_j}$-saturated by the set $I^\ast_{j,0}$. Hence, we can decompose the set $I^{\ast}_{1,m_1}\times \dots \times I^{\ast}_{n,m_n}$ into pairwise non-overlapping  intervals $\Delta_1\times\dots \times \Delta_n$ each of which is a translate of $R_1\times \dots \times R_n$ and 
$$
\frac{|(\Delta_1\times\dots \times\Delta_n)\cap E|}{|\Delta_1\times\dots \times \Delta_n|}=\frac{|(\Delta_1\times\dots \times\Delta_n)\cap (I^\ast_{1,0}\times\dots\times I^\ast_{n,0})|}{|\Delta_1\times\dots \times \Delta_n|}= 
$$
$$
\frac{|\Delta_1 \cap I^\ast_{1,0}|}{|\Delta_1|} \dots \frac{|\Delta_n \cap I^\ast_{n,0}|}{|\Delta_n|}=\frac{1}{2^{m_1}} \dots  \frac{1}{2^{m_n}}=\frac{1}{2^k}.
$$
Hence, $\{x \in \mathbb{R}^n : M_{\mathcal{B}} (\chi_{E}) (x)\geq 1/2^k\}\supset   
I^{\ast}_{1,m_1}\times \dots \times I^{\ast}_{n,m_n}$. Consequently, taking into account that  $(m_1,\dots, m_n)$ is arbitrary in $\Omega$, we conclude $(4)$ holds.
\\

For any $j\in \{1,\dots, n\}$ we denote  $H_{j,0}=I^{\ast}_{j,0}$ and $H_{j,m}=I^{\ast}_{j,m}\setminus I^{\ast}_{j,m-1} $ for $m\in \{1,\dots, k\}$. By virtue of Lemma \ref{l2} it is easy to see that:
\\

\noindent i) The sets $H_{1,m_1}\times\dots\times H_{n,m_n}\;\; ((m_1,\dots, m_n)\in \Omega_{n,k})$ are pairwise disjoint;
\\

\noindent ii) For every $j\in\{1,\dots, n\}$ and $m\in \{1,\dots, k\}$ we have that $|H_{j,m}|=|I^{\ast}_{j,m}|/2$. Consequently, for every 
 $(m_1,\dots, m_n)\in \Omega_{n,k}$,
$$
|H_{1,m_1}\times\dots\times H_{n,m_n}|=\frac{1}{2^n}|I^\ast_{1,m_1}\times\dots\times I^\ast_{n,m_n}|;
$$

\noindent iii) For every 
$(m_1,\dots, m_n)\in \Omega_{n,k}$,
$$
|I^\ast_{1,m_1}\times\dots\times I^\ast_{n,m_n}|= 2^k |I^\ast_{1,0}\times\dots\times I^\ast_{n,0}|.
$$

\noindent Hence,
\begin{align}
&\bigg|\bigcup_{(m_1,\dots, m_n)\in \Omega}  
I^{\ast}_{1,m_1}\times \dots \times I^{\ast}_{n,m_n}\bigg| \notag \\
&\;\;\;\geq \bigg|\bigcup_{(m_1,\dots, m_n)\in \Omega}  
H_{1,m_1}\times \dots \times H_{n,m_n}\bigg| \notag
\\&\;\;\;=\sum_{(m_1,\dots, m_n)\in \Omega}|H_{1,m_1}\times \dots \times H_{n,m_n}| \notag
\\&\;\;\;=\textup{card}\,(\Omega) \frac{1}{2^n} |I^\ast_{1,0}\times\dots\times I^\ast_{n,0}| \notag
\\&\;\;\;= \frac{1}{2^n} \textup{card}\,(\Omega)  2^k |E|. \label{e5}
\end{align}

From (\ref{e4}) and (\ref{e5}) we conclude the  lemma holds.
\end{proof}
\medskip

For an  interval $R$  we denote by $R_d$ the smallest interval concentric to $R$  which contains $R$ and has dyadic side lengths.

Let $\mathcal{B}$ be a rare basis in $\mathbb{R}^n$. To $\mathcal{B}$ we can associate its \textit{dyadic skeleton} $\mathcal{B}_d=\{R_d: R\in \mathcal{B}\}$. Note that the maximal operators associated with the bases $\mathcal{B}$ and $\mathcal{B}_d$ possess similar properties. Namely, $M_\mathcal{B} f\leq 2^n M_{\mathcal{B}_d} f$ and on the other hand (see, e.g., \cite{onianichubinidze},  Lemma 2.12), 
\begin{equation}\label{e6}
|\{x \in \mathbb{R}^n : M_{\mathcal{B}_d} f(x)>\alpha\}|\leq C_n |\{x \in \mathbb{R}^n : M_{\mathcal{B}} f (x)>\alpha/4^n\}|. 
\end{equation}

\medskip
 Theorem \ref{t1} follows from  Lemmas  \ref{l1} and \ref{l3} and estimate (\ref{e6}).


\section{A Generalization of Theorem \ref{t1} and  Applications}

Let $k \geq n$ and $\Omega\subset \Omega_{n,k}$. We will  say that a rare basis  $\mathcal{B}$  in $\mathbb{R}^n$ is $\Omega$-\textit{complete}  if there exist increasing  sequences of integers  $(s_{1,m})_{m=0}^k, \dots,$ $(s_{n,m})_{m=0}^k$ such that for every  $n$-tuple  $(m_1,\dots, m_n)$ belonging to $\Omega$ there exists an interval $R_1\times \dots \times R_n \in \mathcal{B}_d$ with $
 |R_1|\in(2^{s_{1, m_1-1}}, 2^{s_{1, m_1}}],\dots,|R_n|\in(2^{s_{n, m_n-1}}, 2^{s_{n, m_n}}]$.

From  Lemma \ref{l3} and estimate (\ref{e6})  we obtain the following result.

\medskip

\begin{thm}\label{t2} \textit{Let $\mathcal{B}$ be a  rare basis  in $\mathbb{R}^n$. Suppose there exist an increasing sequence of natural numbers $k_j\geq n$ and a sequence of  sets $\Omega_j\subset \Omega_{n,k_j}$ with the properties: $\mathcal{B}$ is $\Omega_j$-complete for every $j\in\mathbb{N}$ and $\inf_{j\in \mathbb{N}}\textup{card}(\Omega_j)/k_j^{n-1} >0$. Then 	 $M_\mathcal{B}$  satisfies the sharp weak type $L(1+ \log^{+}L)^{n-1}$  estimate. Moreover, for every $j \in \mathbb{N}$ there exists a bounded set $E_{j}\subset \mathbb{R}^n$ with  positive  measure such that
$$
|\{x \in \mathbb{R}^n : M_{\mathcal{B}} (\chi_{E_{j}})(x) \geq 1/2^{k_j}\}|\geq c_{\mathcal{B}} k_j^{n-1}2^{k_j}|E_{j}|,
$$
where  $c_\mathcal{B}$ is a constant  of the  form $c_n \inf_{j\in \mathbb{N}}\textup{card}(\Omega_j)/k_j^{n-1}$.
}
\end{thm}
\smallskip

Theorem \ref{t2} is an extension of Theorem \ref{t1} since by  Lemma \ref{l1} the density of the spectrum $W_\mathcal{B}$  in the Cartesian  product  of  infinite sets  $S_1, \dots, S_n \subset \mathbb{Z}$ implies $\Omega_{n,k}$-completeness of the basis $\mathcal{B}$  for every $k\geq n$.



\mbox{ }
\\

We now indicate eight applications of Theorems \ref{t1} and \ref{t2} to rare bases.  
\medskip

\textbf{I}.  Let  $S_1,\dots, S_n \subset \mathbb{Z}$ be infinite sets and $\mathcal{B}$ be the  basis consisting of all   $n$-dimensional intervals with side lengths of the form $2^{s_1}, \dots, 2^{s_n}$ where $s_1,\dots, s_n$ belong to the sets  $S_1, \dots, S_n$ respectively.  Taking into account that the spectrum of $\mathcal{B}$ is the product $S_1\times\dots \times S_n$ and applying Theorem \ref{t1}  for $\mathcal{B}$ and $S_1,\dots, S_n$ we obtain the result proved in \cite{stokolos2006}.

\medskip
\textbf{II (Soria Bases)}.  Let $\Gamma\subset \mathbb{Z}$ be an infinite set  and let $\mathcal{B}$ be the  basis of all $3$-dimensional intervals $R_1\times R_2 \times R_3$ such that $|R_1|, |R_2|\in \mathbb{D}$ and $|R_3|=2^{\gamma}/|R_2|$ for some $\gamma \in \Gamma$, where here and in later applications we denote  the set of  dyadic numbers $\{2^s: s\in \mathbb{Z}\}$ by $\mathbb{D}$. It is easy to see that the spectrum $W_\mathcal{B}$ is the  set
$$
\{(w_1,w_2, w_3 ):w_1,w_2\in \mathbb{Z}, w_3\in \Gamma - w_2\}
$$
and  $W_\mathcal{B}$ is dense in $\mathbb{Z}\times \mathbb{Z} \times \Gamma$. Hence taking  $S_1=S_2=\mathbb{Z}$ and $S_3=\Gamma$ by Theorem \ref{t1} we obtain the sharp weak type  $L(1+ \log^{+}L)^{2}$ estimate for the maximal operator $M_{\mathcal{B}}$ associated to the basis $\mathcal{B}$ which was proved in \cite{dhs2022}.

\medskip
\textbf{III (Zygmund Bases).} Let $\Gamma\subset \mathbb{Z}$ be an infinite set  and let $\mathcal{B}$ be the  basis of all $3$-dimensional intervals $R_1\times R_2 \times R_3$ such that $|R_1|, |R_2|\in \mathbb{D}$ and $|R_3|=2^{\gamma}|R_2|$ for some $\gamma \in \Gamma$. It is easy to see that the spectrum $W_\mathcal{B}$ is the  set
$$
\{(w_1,w_2, w_3 ):w_1,w_2\in \mathbb{Z}, w_3\in \Gamma + w_2\}
$$
and  $W_\mathcal{B}$ is dense in $\mathbb{Z}\times \mathbb{Z} \times \Gamma$. Hence taking  $S_1=S_2=\mathbb{Z}$ and $S_3=\Gamma$ by Theorem \ref{t1} we obtain the sharp weak type  $L(1+ \log^{+}L)^{2}$ estimate for the maximal operator $M_{\mathcal{B}}$ associated to the basis $\mathcal{B}$ which was proved in \cite{hs2022pams}.

\medskip
\textbf{IV (C\'ordoba Bases)}. Let $\Gamma\subset \mathbb{Z}$ be an infinite set  and let $\mathcal{B}$ be the  basis of all $3$-dimensional intervals $R_1\times R_2 \times R_3$ such that $|R_1|, |R_2|\in \mathbb{D}$ and $|R_3|=2^{\gamma}|R_1||R_2|$ for some $\gamma \in \Gamma$. It is easy to see that the spectrum $W_\mathcal{B}$ is the  set
$$
\{(w_1,w_2, w_3 ):w_1,w_2\in \mathbb{Z}, w_3\in \Gamma +w_1 + w_2\}.
$$
and  $W_\mathcal{B}$ is dense in $\mathbb{Z}\times \mathbb{Z} \times \Gamma$. Hence, taking  $S_1=S_2=\mathbb{Z}$ and $S_3=\Gamma$ by Theorem \ref{t1} we obtain the sharp weak type  $L(1+ \log^{+}L)^{2}$ estimate for the maximal operator $M_{\mathcal{B}}$ associated to the basis $\mathcal{B}$ which was proved in \cite{hs2022cordoba}.

\medskip

\medskip
\textbf{V}. Suppose  $T_1,\dots, T_{n-1}$  are infinite subsets of $\mathbb{D}$, $\Gamma$ is an infinite subset of $\mathbb{Z}$,  $1\leq p \leq n-1$, and $1\leq j_1<\dots< j_p\leq n-1$.  Let  $\mathcal{B}$ be the basis of all $n$-dimensional intervals $R_1\times \dots \times R_n$ such that $|R_1|\in T_1, \dots, |R_{n-1}|\in T_{n-1}$ and $|R_n|=2^{\gamma}/(|R_{j_1}|\dots |R_{j_p}|)$ for some $\gamma \in \Gamma$.

Let $S_j=\{\log k: k\in T_j\}$ $(j\in \{1,\dots, n-1\})$ and  $S_{n}=\Gamma$. 

It is easy to see that the spectrum $W_\mathcal{B}$ is the  set
$$
\{(w_1,\dots, w_{n-1}, w_n):w_1\in S_1, \dots,  w_{n-1}\in S_{n-1}, w_n\in \Gamma-(w_{j_1}+\dots + w_{j_p})\}
$$
and  $W_\mathcal{B}$ is dense in $S_1\times\dots \times S_{n}$.

Applying Theorem \ref{t1} for $\mathcal{B}$ and $S_1,\dots,S_{n}$ we obtain the sharp weak type $L(1+ \log^{+}L)^{n-1}$  estimate for the maximal operator $M_{\mathcal{B}}$ associated to the basis $\mathcal{B}$.

Under the same conditions we can obtain the sharp weak type $L(1+ \log^{+}L)^{n-1}$   estimate for the maximal operator $M_{\mathcal{B}}$ associated to the basis  $\mathcal{B}$ of all intervals $R_1\times \dots \times R_n$ such that $|R_1|\in T_1, \dots, |R_{n-1}|\in T_{n-1}$ and $|R_n|=2^{\gamma}|R_{j_1}|\dots |R_{j_p}|$ for some $\gamma \in \Gamma$.

Note that the bases considered in this application are multi-dimensional versions of ones from  applications II-IV.

\medskip
\textbf{VI}. The conditions of Theorem \ref{t1} are satisfied by more general bases than ones considered in the applications II-V. In particular,  let  $\theta_k : \mathbb{D}^{n-1}\rightarrow \mathbb{D}$ $(k\in \mathbb{N})$ be functions  satisfying the following conditions:

1) $\inf\limits_{k\in \mathbb{N}}\theta_k(1,\dots,1)=0\;$ or  $\;\sup\limits_{k\in \mathbb{N}}\theta_k(1,\dots,1)=\infty$;

2) for every $(t_1,\dots, t_{n-1})\in\mathbb{D}^{n-1}$
$$
\inf_{k\in \mathbb{N}}\frac{\theta_k(t_1,\dots,t_{n-1})}{\theta_k(1,\dots,1)}>0 \;\;\text{and}\;\;\sup_{k\in \mathbb{N}}\frac{\theta_k(t_1,\dots,t_{n-1})}{\theta_k(1,\dots,1)}<\infty.
$$
Suppose  $T_1,\dots, T_{n-1}$  are infinite subsets of $\mathbb{D}$. Let  $\mathcal{B}$ be the basis of all $n$-dimensional intervals $R_1\times \dots \times R_n$ such that $|R_1|\in T_1, \dots, |R_{n-1}|\in T_{n-1}$ and $|R_n|=\theta_k(|R_1|,\dots,|R_{n-1}|)$ for some $k\in \mathbb {N}$.

Let $S_j=\{\log k: k\in T_j\}\; (j\in \{1,\dots, n-1\})$ and $S_{n}=\{\log \theta_k(1,\dots,1): k\in \mathbb{N}\}$.

It is easy to see that the spectrum $W_\mathcal{B}$ is the following set
$$
\{(w_1,\dots, w_{n-1}, \log \theta_k(2^{w_1},\dots, 2^{w_{n-1}})): w_1\in S_1, \dots, w_{n-1}\in S_{n-1}, k\in\mathbb{N}\}
$$
and  $W_\mathcal{B}$ is dense in $S_1\times\dots \times S_{n}$.

Applying Theorem \ref{t1} for $\mathcal{B}$ and  $S_1,\dots,S_{n}$ we obtain the sharp weak type $L(1+ \log^{+}L)^{n-1}$   estimate for the maximal operator $M_{\mathcal{B}}$ associated to the basis $\mathcal{B}$.

\medskip
\textbf{VII}. Following \cite{stok2008} let us say that a rare basis $\mathcal{B}$ in $\mathbb{R}^2$  satisfies the (is)-\textit{property}  if for every $k\in \mathbb{N}$ there exist intervals $R_0,\dots, R_k\in \mathcal{B}$ of the type $[0,2^p]\times [0,2^q]$ $(p,q\in \mathbb{Z})$ such that:

1) For every $i,j\in\{0,\dots, k\}$ with $i\neq j$ the intervals $R_i$ and $R_j$ are incomparable,  i.e., there does not exist translation placing one of them inside the other; 

2) For every $i,j\in\{0,\dots, k\}$  the interval $R_i\cap R_j$ belongs to  $\mathcal{B}$. 

Let $\mathcal{B}_1$ be a basis of two-dimensional intervals with the (is)-property and $\mathcal{B}_2$ be a basis consisting of one-dimensional intervals with  lengths belonging to an infinite set of dyadic numbers. By $\mathcal{B}_1\times \mathcal{B}_2$ denote their  product, i.e., the basis which consists of three-dimensional intervals of the type 
$J_1\times J_2$ where $J_1\in \mathcal{B}_1$ and $J_2\in \mathcal{B}_2$. 

Let  $k\geq 3$. We can find intervals $R_0,\dots, R_k\in \mathcal{B}_1$ of the type $[0,2^p]\times [0,2^q]$ $(p,q\in \mathbb{Z})$  with the properties 1) and 2) from the definition of the (is)-property.  We can assume that
$$
R_0=[0,2^{p_0}]\times [0,2^{q_{k}}],\dots,R_i=[0,2^{p_i}]\times [0,2^{q_{k-i}}], \dots, R_k=[0,2^{p_k}]\times [0,2^{q_{0}}],
$$  
where $p_0<\dots <p_k $ and $q_0<\dots < q_k $. Let $t_0<\dots <t_k $ be integers such that $[0,2^{t_0}], \dots, [0,2^{t_k}]\in \mathcal{B}_2$. Set $s_{1,0}=p_0, \dots, s_{1,k}=p_k,$ $s_{2,0}=q_0, \dots, s_{2,k}=q_k$,  and $s_{3,0}=t_0, \dots, s_{3,k}=t_k$.     Then for every triple $(m_1, m_2, m_3)$ belonging to $V_{3,k}$ it is easy to see that
$$
[0,2^{s_{1,m_1}}]\times [0,2^{s_{2,m_2}}]\times [0,2^{s_{3,m_3}}]=(R_{m_1}\cap R_{k-m_2})\times [0,2^{t_{m_3}}] \in  \mathcal{B}_1\times \mathcal{B}_2.
$$ 
Hence, $\mathcal{B}_1\times \mathcal{B}_2$ is $\Omega_{3,k}$-complete for every $k\geq 3$.

Applying Theorem \ref{t2}  we obtain the sharp weak type $L(1+ \log^{+}L)^{2}$  estimate for the maximal operator $M_{\mathcal{B}}$ associated to the basis $\mathcal{B}$ where  $\mathcal{B} = \mathcal{B}_1\times \mathcal{B}_2$.   For the case of  $\mathcal{B}_2$ being the basis of all intervals with dyadic lengths the estimate was obtained in \cite{stok2008}.

\medskip
\textbf{VIII}. In  Section 4.2 of \cite{dm2017} a certain  class  $\Lambda_n$ of rare bases in $\mathbb{R}^n$ is  considered such that every basis $\mathcal{B}$ from $\Lambda_n$ has the following property (see Remark 9 in \cite{dm2017}):  For every $k \in \mathbb{N}$ there exist  intervals $R_0,R_1,\dots, R_k$ such that
$$
R_0=[0,2^{s_{1,0}}]\times[0,2^{s_{2,0}}]\times \dots\times [0,2^{s_{n,0}}],
$$
$$
R_1=[0,2^{s_{1,1}}]\times[0,2^{s_{2,1}}]\times\dots\times [0,2^{s_{n,1}}],
$$
$$\hspace{-2in}\vdots$$
$$
R_k=[0,2^{s_{1,k}}]\times[0,2^{s_{2,k}}]\times\dots\times [0,2^{s_{n,k}}],
$$
where $(s_{1,m})_{m=0}^k, \dots, (s_{n,m})_{m=0}^k$ are increasing sequences of integers and for every $n$-tuple of  integers $(m_1,\dots,m_n)$ with   $k\geq m_1\geq m_2 \geq \dots \geq m_n\geq 0$ the interval
$$
R=[0,2^{s_{1,m_1}}]\times[0,2^{s_{2,m_2}}]\times \dots\times [0,2^{s_{n,m_n}}]
$$
belongs to the basis $\mathcal{B}$.

 Suppose  $\mathcal{B}_1$  is a basis  from the class  $\Lambda_{n}$ and $\mathcal{B}_2$ is a basis consisting of one-dimensional intervals with  lengths belonging to an infinite set of dyadic numbers.  Taking into account the above given property of bases from the class $\Lambda_n$  we have that the product basis $\mathcal{B}_1 \times \mathcal{B}_2$ is $\Omega_k$-complete for every $k\geq n+1$ where  $\Omega_k$ is the set of all $(n+1)$-tuples $(m_1,\dots,m_n, m_{n+1})\in \Omega_{n+1,k}$ with $k\geq m_1\geq m_2 \geq \dots \geq m_{n}\geq 1$. On the other hand, it is easy to see that  $\textup{card}\, \Omega_k \geq c_n k^{n}$ $(k\geq n+1)$.
 
  Applying Theorem \ref{t2} we obtain the sharp weak type $L(1+ \log^{+}L)^{n}$   estimate for the maximal operator $M_{\mathcal{B}}$ associated to the basis $\mathcal{B}$ where $\mathcal{B} = \mathcal{B}_1\times \mathcal{B}_2$.

\end{document}